\documentclass[10pt, a4paper,reqno]{amsart}
\usepackage{amsfonts}
\usepackage{amssymb, amsmath}
\usepackage{bbm}
\numberwithin{equation}{section}
\usepackage{url}

\newtheorem{theorem}{Theorem}[section]
\newtheorem{proposition}[theorem]{Proposition}
\newtheorem{corollary}[theorem]{Corollary}

\theoremstyle{definition}

\newtheorem{example}[theorem]{\textbf{Example}}
\newtheorem{remark}[theorem]{\textbf{Remark}}

\title[Complete intersection]{Complete intersection monomial curves and the Cohen-Macaulayness of their tangent cones}

\author[A. Katsabekis]{Anargyros Katsabekis}
\address {Department of Mathematics, Bilkent University, 06800 Ankara, Turkey} \email{katsampekis@bilkent.edu.tr}


\begin{document}

\begin{abstract} Let $C({\bf n})$ be a complete intersection monomial curve in the 4-dimensional affine space. In this paper we study the complete intersection property of the monomial curve $C({\bf n}+w{\bf v})$, where $w>0$ is an integer and ${\bf v} \in \mathbb{N}^{4}$. Also we investigate the Cohen-Macaulayness of the tangent cone of $C({\bf n}+w{\bf v})$.

\end{abstract}

\maketitle

\section{Introduction}

Let ${\bf n}=(n_{1},n_{2},\ldots,n_{d})$ be a sequence of positive integers with ${\rm gcd}(n_{1},\ldots,n_{d})=1$. Consider the polynomial ring $K[x_{1},\ldots,x_{d}]$ in $d$ variables over a field $K$. We shall denote by ${\bf x}^{\bf u}$ the monomial $x_{1}^{u_1} \cdots x_{d}^{u_d}$ of $K[x_{1},\ldots,x_{d}]$, with ${\bf u}=(u_{1},\ldots,u_{d}) \in \mathbb{N}^{d}$ where $\mathbb{N}$ stands for the set of non-negative integers. The toric ideal $I({\bf n})$ is the kernel of the $K$-algebra homomorphism $\phi:K[x_{1},\ldots,x_{d}] \rightarrow K[t]$ given by $$\phi(x_{i})=t^{n_i} \ \ \textrm{for all} \ \ 1 \leq i \leq d.$$ Then $I({\bf n})$ is the defining ideal of the monomial curve $C({\bf n})$ given by the parametrization $x_{1}=t^{n_1},\ldots,x_{d}=t^{n_d}$. The ideal $I({\bf n})$ is generated by all the binomials ${\bf x}^{{\bf u}}-{\bf x}^{{\bf v}}$, where ${\bf u}-{\bf v}$ runs over all vectors in the lattice ${\rm ker}_{\mathbb{Z}}(n_{1},\ldots,n_{d})$ see for example, \cite[Lemma 4.1]{Sturmfels95}. The height of $I({\bf n})$ is $d-1$ and also equals the rank of ${\rm ker}_{\mathbb{Z}}(n_{1},\ldots,n_{d})$ (see \cite{Sturmfels95}). Given a polynomial $f \in I({\bf n})$, we let $f_{*}$ be the homogeneous summand of $f$ of least degree. We shall denote by $I({\bf n})_{*}$ the ideal in $K[x_{1},\ldots,x_{d}]$ generated by the polynomials $f_{*}$ for $f \in I({\bf n})$.

Deciding whether the associated graded ring of the local ring $K[[t^{n_{1}},\ldots,t^{n_{d}}]]$ is Cohen-Macaulay constitutes an important problem studied by many authors, see for instance \cite{AMS}, \cite{Ga}, \cite{Shi}. The importance of this problem stems partially from the fact that if the associated graded ring is Cohen-Macaulay, then the Hilbert function of $K[[t^{n_{1}},\ldots,t^{n_{d}}]]$ is non-decreasing. Since the associated graded ring of $K[[t^{n_{1}},\ldots,t^{n_{d}}]]$ is isomorphic to the ring $K[x_{1},\ldots,x_{d}]/I({\bf n})_{*}$, the Cohen-Macaulayness of the associated graded ring can be studied as the Cohen-Macaulayness of the ring $K[x_{1},\ldots,x_{d}]/I({\bf n})_{*}$. Recall that $I({\bf n})_{*}$ is the defining ideal of the tangent cone of $C({\bf n})$ at $0$.

The case that $K[[t^{n_{1}},\ldots,t^{n_{d}}]]$ is Gorenstein has been particularly studied. This is partly due to the M. Rossi's problem \cite{Rossi} asking whether the Hilbert function of a Gorenstein local ring of dimension one is non-decreasing. Recently, A. Oneto, F. Strazzanti and G. Tamone \cite{OST} found many families of monomial curves giving negative answer to the above problem. However M. Rossi's problem is still open for a Gorenstein local ring $K[[t^{n_{1}},\ldots,t^{n_{4}}]]$. It is worth to note that, for a complete intersection monomial curve $C({\bf n})$ in the 4-dimensional affine space (i.e. the ideal $I({\bf n})$ is a complete intersection), we have, from \cite[Theorem 3.1]{Shi}, that if the minimal number of generators for $I({\bf n})_{*}$ is either three or four, then $C({\bf n})$ has Cohen-Macaulay tangent cone at the origin. The converse is not true in general, see \cite[Proposition 3.14]{Shi}.

In recent years there has been a surge of interest in studying properties of the monomial curve $C({\bf n}+w{\bf v})$, where $w>0$ is an integer and ${\bf v} \in \mathbb{N}^{d}$, see for instance \cite{CGHN}, \cite{Gi} and \cite{Vu}. This is particularly true for the case that ${\bf v}=(1,\ldots,1)$. In fact, J. Herzog and H. Srinivasan conjectured that if $n_{1}<n_{2}<\cdots <n_{d}$ are positive numbers, then the Betti numbers of $I({\bf n}+w{\bf v})$ are eventually periodic in $w$ with period $n_{d}-n_{1}$. The conjecture was proved by T. Vu \cite{Vu}. More precisely, he showed that there exists a positive integer $N$ such that, for all $w>N$, the Betti numbers of $I({\bf n}+w{\bf v})$ are periodic in $w$ with period $n_{d}-n_{1}$. The bound $N$ depends on the Castelnuovo-Mumford regularity of the ideal generated by the homogeneous elements in $I({\bf n})$. For $w>(n_{d}-n_{1})^{2}-n_{1}$ the minimal number of generators for $I({\bf n}+w(1,\ldots,1))$ is periodic in $w$ with period $n_{d}-n_{1}$ (see \cite{CGHN}). Furthermore, for every $w>(n_{d}-n_{1})^{2}-n_{1}$ the monomial curve $C({\bf n}+w(1,\ldots,1))$ has Cohen-Macaulay tangent cone at the origin, see  \cite{Stamate}. The next example provides a monomial curve $C({\bf n}+w(1,\ldots,1))$ which is not a complete intersection for every $w>0$.

\begin{example} {\rm Let ${\bf n}=(15,25,24,16)$, then $I({\bf n})$ is a complete intersection on the binomials $x_{1}^{5}-x_{2}^{3}$, $x_{3}^{2}-x_{4}^{3}$ and $x_{1}x_{2}-x_{3}x_{4}$. Consider the vector ${\bf v}=(1,1,1,1)$. For every $w>85$ the minimal number of generators for $I({\bf n}+w{\bf v})$ is either $18$, $19$ or $20$. Using CoCoA (\cite{Coc}) we find that for every $0<w \leq 85$ the minimal number of generators for $I({\bf n}+w{\bf v})$ is greater than or equal to $4$. Thus for every $w>0$ the ideal $I({\bf n}+w{\bf v})$ is not a complete intersection.}
\end{example}

Given a complete intersection monomial curve $C({\bf n})$ in the 4-dimensional affine space, we study (see Theorems \ref{Basic-complete}, \ref{Basic-complete3}) when $C({\bf n}+w{\bf v})$ is a complete intersection. We also construct (see Theorems \ref{Basic1}, \ref{TangentAll}, \ref{theorcomp}) families of complete intersection monomial curves $C({\bf n}+w{\bf v})$ with Cohen-Macaulay tangent cone at the origin.

Let $a_{i}$ be the least positive integer such that $a_{i}n_{i} \in \sum_{j \neq i} \mathbb{N}n_{j}$. To study the complete intersection property of $C({\bf n}+w{\bf v})$ we use the fact that after permuting variables, if necessary, there exists (see \cite[Proposition 3.2]{Shi} and also Theorems 3.6 and 3.10 in \cite{KaOj}) a minimal system
of binomial generators $S$ of $I({\bf n})$ of the following form:  \begin{enumerate} \item[(A)] $S=\{x_{1}^{a_1}-x_{2}^{a_2}, x_{3}^{a_{3}}-x_{4}^{a_4}, x_{1}^{u_{1}}x_{2}^{u_{2}}-x_{3}^{u_{3}}x_{4}^{u_{4}}\}$. \item[(B)] $S=\{x_{1}^{a_1}-x_{2}^{a_2}, x_{3}^{a_{3}}-x_{1}^{u_1}x_{2}^{u_2}, x_{4}^{a_{4}}-x_{1}^{v_1}x_{2}^{v_2}x_{3}^{v_3}\}$. 
\end{enumerate}

In section 2 we focus on case (A). We prove that the monomial curve $C({\bf n})$ has Cohen-Macaulay tangent cone at the origin if and only if the minimal number of generators for $I({\bf n})_{*}$ is either three or four. Also we explicitly construct vectors ${\bf v}_{i}$, $1 \leq i \leq 22$, such that for every $w>0$, the ideal $I({\bf n}+w{\bf v}_{i})$ is a complete intersection whenever the entries of ${\bf n}+w{\bf v}_{i}$ are relatively prime. We show that if $C({\bf n})$ has Cohen-Macaulay tangent cone at the origin, then for every $w>0$ the monomial curve $C({\bf n}+w{\bf v}_{1})$ has Cohen-Macaulay tangent cone at the origin whenever the entries of ${\bf n}+w{\bf v}_{1}$ are relatively prime. Additionally we show that there exists a non-negative integer $w_{0}$ such that for all $w \geq w_{0}$, the monomial curves $C({\bf n}+w{\bf v}_{9})$ and $C({\bf n}+w{\bf v}_{13})$ have Cohen-Macaulay tangent cones at the origin whenever the entries of the corresponding sequence (${\bf n}+w{\bf v}_{9}$ for the first family and ${\bf n}+w{\bf v}_{13}$ for the second) are relatively prime. Finally we provide an infinite family of complete intersection monomial curves $C_{m}({\bf n}+w{\bf v}_{1})$ with corresponding local rings having non-decreasing Hilbert functions, although their tangent cones are not Cohen-Macaulay, thus giving a positive partial answer to M. Rossi's problem. 

In section 3 we study the case (B). We construct vectors ${\bf b}_{i}$, $1 \leq i \leq 22$, such that for every $w>0$, the ideal $I({\bf n}+w{\bf b}_{i})$ is a complete intersection whenever the entries of ${\bf n}+w{\bf b}_{i}$ are relatively prime. Furthermore we show that there exists a non-negative integer $w_{1}$ such that for all $w \geq w_{1}$, the ideal $I({\bf n}+w{\bf b}_{22})_{*}$ is a complete intersection whenever the entries of ${\bf n}+w{\bf b}_{22}$ are relatively prime.

\section{The case (A)} 

In this section we assume that after permuting variables, if necessary, $S=\{x_{1}^{a_1}-x_{2}^{a_2}, x_{3}^{a_{3}}-x_{4}^{a_4}, x_{1}^{u_{1}}x_{2}^{u_{2}}-x_{3}^{u_{3}}x_{4}^{u_{4}}\}$ is a minimal generating set of $I({\bf n})$. First we will show that the converse of \cite[Theorem 3.1]{Shi} is also true in this case.

Let $n_{1}={\rm min}\{n_{1},\ldots,n_{4}\}$ and also $a_{3}<a_{4}$. By \cite[Theorem 7]{Ga} a monomial curve $C({\bf n})$ has Cohen-Macaulay tangent cone if and only if $x_{1}$ is not a zero divisor in the ring $K[x_{1},\ldots,x_{4}]/I({\bf n})_{*}$. Hence if $C({\bf n})$ has Cohen-Macaulay tangent cone at the origin, then $I({\bf n})_{*}: \langle x_{1}\rangle=I({\bf n})_{*}$. Without loss of generality we can assume that $u_{2} \leq a_{2}$. In case that $u_{2}>a_{2}$ we can write $u_{2}=ga_{2}+h$, where $0 \leq h<a_{2}$. Then we can replace the binomial $x_{1}^{u_{1}}x_{2}^{u_{2}}-x_{3}^{u_{3}}x_{4}^{u_{4}}$ in $S$ with the binomial $x_{1}^{u_{1}+ga_{1}}x_{2}^{h}-x_{3}^{u_{3}}x_{4}^{u_{4}}$. Without loss of generality we can also assume that $u_{3} \leq a_{3}$.

\begin{theorem} \label{CM-tangent} Suppose that $u_{3}>0$ and $u_{4}>0$. Then $C({\bf n})$ has Cohen-Macaulay tangent cone at the origin if and only if the ideal $I({\bf n})_{*}$ is either a complete intersection or an almost complete intersection.

\end{theorem}

\noindent \textbf{Proof.} $(\Longleftarrow)$ If the minimal number of generators of $I({\bf n})_{*}$ is either three or four, then $C({\bf n})$ has Cohen-Macaulay tangent cone at the origin.

$(\Longrightarrow)$ Let $f_{1}=x_{1}^{a_1}-x_{2}^{a_2}$, $f_{2}=x_{3}^{a_{3}}-x_{4}^{a_4}$, $f_{3}=x_{1}^{u_{1}}x_{2}^{u_{2}}-x_{3}^{u_{3}}x_{4}^{u_{4}}$. We distinguish the following cases \begin{enumerate} \item $u_{2}<a_{2}$. Note that $x_{4}^{a_{4}+u_{4}}-x_{1}^{u_{1}}x_{2}^{u_{2}}x_{3}^{a_{3}-u_{3}} \in I({\bf n})$. We will show that $a_{4}+u_{4} \leq u_{1}+u_{2}+a_{3}-u_{3}$. Suppose that $u_{1}+u_{2}+a_{3}-u_{3}<a_{4}+u_{4}$, then $x_{2}^{u_{2}}x_{3}^{a_{3}-u_{3}} \in I({\bf n})_{*}: \langle x_{1}\rangle$ and therefore $x_{2}^{u_{2}}x_{3}^{a_{3}-u_{3}} \in I({\bf n})_{*}$. Since $\{f_{1},f_{2},f_{3}\}$ is a generating set of $I({\bf n})$, the monomial $x_{2}^{u_{2}}x_{3}^{a_{3}-u_{3}}$ is divided by at least one of the monomials $x_{2}^{a_2}$ and $x_{3}^{a_{3}}$. But $u_{2}<a_{2}$ and $a_{3}-u_{3}<a_{3}$, so $a_{4}+u_{4} \leq u_{1}+u_{2}+a_{3}-u_{3}$. Let $G=\{f_{1}, f_{2}, f_{3}, f_{4}=x_{4}^{a_{4}+u_{4}}-x_{1}^{u_{1}}x_{2}^{u_{2}}x_{3}^{a_{3}-u_{3}}\}$. We will prove that $G$ is a standard basis for $I({\bf n})$ with respect to the negative degree reverse lexicographical order with $x_{3}>x_{4}>x_{2}>x_{1}$. Note that $u_{3}+u_{4}<u_{1}+u_{2}$, since $u_{3}+u_{4} \leq u_{1}+u_{2}+a_{3}-a_{4}$ and also $a_{3}-a_{4}<0$. Thus ${\rm LM}(f_{3})=x_{3}^{u_{3}}x_{4}^{u_{4}}$. Furthermore ${\rm LM}(f_{1})=x_{2}^{a_2}$, ${\rm LM}(f_{2})=x_{3}^{a_{3}}$ and ${\rm LM}(f_{4})=x_{4}^{a_{4}+u_{4}}$. Therefore ${\rm NF}({\rm spoly}(f_{i},f_{j})|G) = 0$ as ${\rm LM}(f_{i})$ and ${\rm LM}(f_j)$ are relatively prime, for $(i,j) \in \{(1,2),(1,3),(1,4),(2,4)\}$. We compute ${\rm spoly}(f_{2},f_{3})=-f_{4}$, so ${\rm NF}({\rm spoly}(f_{2},f_{3})|G)=0$. Next we compute ${\rm spoly}(f_{3},f_{4})=x_{1}^{u_{1}}x_{2}^{u_{2}}x_{3}^{a_{3}}-x_{1}^{u_{1}}x_{2}^{u_{2}}x_{4}^{a_{4}}$. Then ${\rm LM}({\rm spoly}(f_{3},f_{4}))=x_{1}^{u_{1}}x_{2}^{u_{2}}x_{3}^{a_{3}}$ and only ${\rm LM}(f_{2})$ divides ${\rm LM}({\rm spoly}(f_{3},f_{4}))$. Also ${\rm ecart}({\rm spoly}(f_{3},f_{4}))=a_{4}-a_{3}={\rm ecart}(f_{2})$. Then ${\rm spoly}(f_{2},{\rm spoly}(f_{3},f_{4}))=0$ and ${\rm NF}({\rm spoly}(f_{3},f_{4})|G)=0$. By \cite[Lemma 5.5.11]{GP} $I({\bf n})_*$ is generated by the least homogeneous summands of the elements in the standard basis $G$. Thus the minimal number of generators for $I({\bf n})_*$ is least than or equal to 4.

\item $u_{2}=a_{2}$. Note that $x_{4}^{a_{4}+u_{4}}-x_{1}^{u_{1}+a_{1}}x_{3}^{a_{3}-u_{3}} \in I({\bf n})$. We will show that $a_{4}+u_{4} \leq u_{1}+a_{1}+a_{3}-u_{3}$. Clearly the above inequality is true when $u_{3}=a_{3}$. Suppose that $u_{3}<a_{3}$ and $u_{1}+a_{1}+a_{3}-u_{3}<a_{4}+u_{4}$, then $x_{3}^{a_{3}-u_{3}} \in I({\bf n})_{*}: \langle x_{1} \rangle$ and therefore $x_{3}^{a_{3}-u_{3}} \in I({\bf n})_{*}$. Thus $x_{3}^{a_{3}-u_{3}}$ is divided by $x_{3}^{a_{3}}$, a contradiction. Consequently $a_{4}+u_{4} \leq u_{1}+a_{1}+a_{3}-u_{3}$. We will prove that $H=\{f_{1},f_{2},f_{5}=x_{1}^{u_{1}+a_{1}}-x_{3}^{u_{3}}x_{4}^{u_{4}},f_{6}=x_{4}^{a_{4}+u_{4}}-x_{1}^{u_{1}+a_{1}}x_{3}^{a_{3}-u_{3}}\}$ is a standard basis for $I({\bf n})$ with respect to the negative degree reverse lexicographical order with $x_{3}>x_{4}>x_{2}>x_{1}$. Here ${\rm LM}(f_{1})=x_{2}^{a_{2}}$, ${\rm LM}(f_{2})=x_{3}^{a_{3}}$, ${\rm LM}(f_{5})=x_{3}^{u_{3}}x_{4}^{u_{4}}$ and ${\rm LM}(f_{6})=x_{4}^{u_{4}+a_{4}}$. Therefore ${\rm NF}({\rm spoly}(f_{i},f_{j})|H)=0$ as ${\rm LM}(f_{i})$ and ${\rm LM}(f_{j})$ are relatively prime, for $(i,j) \in \{(1,2),(1,5),(1,6),(2,6)\}$. We compute ${\rm spoly}(f_{2},f_{5})=-f_{6}$, therefore ${\rm NF}({\rm spoly}(f_{2},f_{5})|H)=0$. Furthermore ${\rm spoly}(f_{5},f_{6})=x_{1}^{u_{1}+a_{1}}x_{3}^{a_{3}}-x_{1}^{u_{1}+a_{1}}x_{4}^{a_{4}}$ and also ${\rm LM}({\rm spoly}(f_{5},f_{6}))=x_{1}^{u_{1}+a_{1}}x_{3}^{a_{3}}$. Only ${\rm LM}(f_{2})$ divides ${\rm LM}({\rm spoly}(f_{5},f_{6}))$ and ${\rm ecart}({\rm spoly}(f_{5},f_{6}))=a_{4}-a_{3}={\rm ecart}(f_{2})$. Then ${\rm spoly}(f_{2},{\rm spoly}(f_{5},f_{6}))=0$ and therefore ${\rm NF}({\rm spoly}(f_{5},f_{6})|H)=0$. By \cite[Lemma 5.5.11]{GP} $I({\bf n})_*$ is generated by the least homogeneous summands of the elements in the standard basis $H$. Thus the minimal number of generators for $I({\bf n})_*$ is least than or equal to 4. \hfill $\square$

\end{enumerate}

\begin{corollary} \label{Arithmetical} Suppose that $u_{3}>0$ and $u_{4}>0$. \begin{enumerate} \item Assume that $u_{2}<a_{2}$. Then $C({\bf n})$ has Cohen-Macaulay tangent cone at the origin if and only if $a_{4}+u_{4} \leq u_{1}+u_{2}+a_{3}-u_{3}$. \item Assume that $u_{2}=a_{2}$. Then $C({\bf n})$ has Cohen-Macaulay tangent cone at the origin if and only if $a_{4}+u_{4} \leq u_{1}+a_{1}+a_{3}-u_{3}$.
\end{enumerate}

\end{corollary}

\begin{theorem} \label{CM-tangent1} Suppose that either $u_{3}=0$ or $u_{4}=0$. Then $C({\bf n})$ has Cohen-Macaulay tangent cone at the origin if and only if the ideal $I({\bf n})_{*}$ is a complete intersection.

\end{theorem}

\noindent \textbf{Proof.} It is enough to show that if $C({\bf n})$ has Cohen-Macaulay tangent cone at the origin, then the ideal $I({\bf n})_{*}$ is a complete intersection. Suppose first that $u_{3}=0$. Then $\{f_{1}=x_{1}^{a_1}-x_{2}^{a_2}, f_{2}=x_{3}^{a_{3}}-x_{4}^{a_4}, f_{3}=x_{4}^{u_4}-x_{1}^{u_{1}}x_{2}^{u_{2}}\}$ is a minimal generating set of $I({\bf n})$. If $u_{2}=a_{2}$, then $\{f_{1}, f_{2}, x_{4}^{u_4}-x_{1}^{u_{1}+a_{1}}\}$ is a standard basis for $I({\bf n})$ with respect to the negative degree reverse lexicographical order with $x_{3}>x_{4}>x_{2}>x_{1}$. By \cite[Lemma 5.5.11]{GP} $I({\bf n})_*$ is a complete intersection. Assume that $u_{2}<a_{2}$. We will show that $u_{4} \leq u_{1}+u_{2}$. Suppose that $u_{4}>u_{1}+u_{2}$, then $x_{2}^{u_{2}}\in I({\bf n})_{*}: \langle x_{1}\rangle$ and therefore $x_{2}^{u_{2}} \in I({\bf n})_{*}$. Thus $x_{2}^{u_{2}}$ is divided by $x_{2}^{a_{2}}$, a contradiction. Then $\{f_{1}, f_{2}, f_{3}\}$ is a standard basis for $I({\bf n})$ with respect to the negative degree reverse lexicographical order with $x_{3}>x_{4}>x_{2}>x_{1}$. Note that ${\rm LM}(f_{1})=x_{2}^{a_2}$, ${\rm LM}(f_{2})=x_{3}^{a_{3}}$ and ${\rm LM}(f_{3})=x_{4}^{u_{4}}$. By \cite[Lemma 5.5.11]{GP} $I({\bf n})_*$ is a complete intersection. Suppose now that $u_{4}=0$, so necessarily $u_{3}=a_{3}$. Then $\{f_{1}, f_{2}, f_{4}=x_{4}^{a_4}-x_{1}^{u_{1}}x_{2}^{u_{2}}\}$ is a minimal generating set of $I({\bf n})$. If $u_{2}=a_{2}$, then $\{f_{1},f_{2},x_{4}^{a_4}-x_{1}^{a_{1}+u_{1}}\}$ is a standard basis for $I({\bf n})$ with respect to the negative degree reverse lexicographical order with $x_{3}>x_{4}>x_{2}>x_{1}$. Thus, from \cite[Lemma 5.5.11]{GP}, $I({\bf n})_*$ is a complete intersection. Assume that $u_{2}<a_{2}$, then $a_{4} \leq u_{1}+u_{2}$ and also $\{f_{1},f_{2},f_{4}\}$ is a standard basis for $I({\bf n})$ with respect to the negative degree reverse lexicographical order with $x_{3}>x_{4}>x_{2}>x_{1}$. From \cite[Lemma 5.5.11]{GP} we deduce that $I({\bf n})_*$ is a complete intersection. \hfill $\square$

\begin{remark} {\rm In case (B) the minimal number of generators of $I({\bf n})_{*}$ can be arbitrarily large even if the tangent cone of $C({\bf n})$ is Cohen-Macaulay, see \cite[Proposition 3.14]{Shi}.}

\end{remark}

Given a complete intersection monomial curve $C({\bf n})$, we next study the complete intersection property of $C({\bf n}+w{\bf v})$. Let $M$ be a non-zero $r \times s$ integer matrix, then there exist an $r \times r$ invertible integer matrix $U$ and an $s \times s$ invertible integer matrix $V$ such that $UMV={\rm diag}(\delta_{1},\ldots,\delta_{m},0,\ldots,0)$ is the diagonal matrix, where $\delta_{j}$ for all $j=1,2,\ldots,m$ are positive integers such that $\delta_{i}|\delta_{i+1}$, $1 \leq i \leq m-1$, and $m$ is the rank of $M$. The  elements $\delta_{1},\ldots,\delta_{m}$ are the invariant factors of $M$. By \cite[Theorem 3.9]{Jac} the product $\delta_{1} \delta_{2} \cdots \delta_{m}$ equals the greatest common divisor of all non-zero $m \times m$ minors of $M$.

The following proposition will be useful in the proof of Theorem \ref{Basic-complete}.

\begin{proposition} \label{Complete1} Let $B=\{f_{1}=x_{1}^{b_1}-x_{2}^{b_2}, f_{2}=x_{3}^{b_3}-x_{4}^{b_4}, f_{3}=x_{1}^{v_{1}}x_{2}^{v_{2}}-x_{3}^{v_{3}}x_{4}^{v_{4}}\}$ be a set of binomials in $K[x_{1},\ldots,x_{4}]$, where $b_{i} \geq 1$ for all $1 \leq i \leq 4$, at least one of $v_{1}$, $v_{2}$ is non-zero and at least one of $v_{3}$, $v_{4}$ is non-zero. Let $n_{1}=b_{2}(b_{3}v_{4}+v_{3}b_{4})$, $n_{2}=b_{1}(b_{3}v_{4}+v_{3}b_{4})$, $n_{3}=b_{4}(b_{1}v_{2}+v_{1}b_{2})$, $n_{4}=b_{3}(b_{1}v_{2}+v_{1}b_{2})$. If ${\rm gcd}(n_{1},\ldots,n_{4})=1$, then $I({\bf n})$ is a complete intersection ideal generated by the binomials $f_{1}$, $f_{2}$ and $f_{3}$.

\end{proposition}

\noindent \textbf{Proof.} Consider the vectors ${\bf d}_{1}=(b_{1},-b_{2},0,0)$, ${\bf d}_{2}=(0,0,b_{3},-b_{4})$ and ${\bf d}_{3}=(v_{1},v_{2},-v_{3},-v_{4})$. Clearly ${\bf d}_{i} \in {\rm ker}_{\mathbb{Z}}(n_{1},\ldots,n_{4})$ for $1 \leq i \leq 3$, so the lattice $L=\sum_{i=1}^{3} \mathbb{Z} {\bf d}_{i}$ is a subset of ${\rm ker}_{\mathbb{Z}}(n_{1},\ldots,n_{4})$. Consider the matrix $$M=\begin{pmatrix}
 b_{1} & 0 & v_{1} \\
 -b_{2} & 0 & v_{2} \\
 0 & b_{3} & -v_{3} \\
 0 & -b_{4} & -v_{4}
\end{pmatrix}.$$ It is not hard to show that the rank of $M$ equals $3$. We will prove that $L$ is saturated, namely the invariant factors $\delta_{1}$, $\delta_{2}$ and $\delta_{3}$ of $M$ are all equal to $1$. The greatest common divisor of all non-zero $3 \times 3$ minors of $M$ equals the greatest common divisor of the integers $n_{1}$, $n_{2}$, $n_{3}$ and $n_{4}$. But ${\rm gcd}(n_{1},\ldots,n_{4})=1$, so $\delta_{1} \delta_{2} \delta_{3}=1$ and therefore $\delta_{1}=\delta_{2}=\delta_{3}=1$. Note that the rank of the lattice ${\rm ker}_{\mathbb{Z}}(n_{1},\ldots,n_{4})$ is 3 and also equals the rank of $L$. By \cite[Lemma 8.2.5]{Vil} we have that $L={\rm ker}_{\mathbb{Z}}(n_{1},\ldots,n_{4})$. Now the transpose $M^{t}$ of $M$ is mixed dominating. Recall that a matrix $P$ is mixed dominating if every row of $P$ has a positive and negative entry and $P$ contains no square submatrix with this property. By \cite[Theorem 2.9]{FS} $I({\bf n})$ is a complete intersection on the binomials $f_{1}$, $f_{2}$ and $f_{3}$. \hfill $\square$

\begin{theorem} \label{Basic-complete} Let $I({\bf n})$ be a complete intersection ideal generated by the binomials $f_{1}=x_{1}^{a_1}-x_{2}^{a_2}$, $f_{2}=x_{3}^{a_3}-x_{4}^{a_4}$ and $f_{3}=x_{1}^{u_1}x_{2}^{u_2}-x_{3}^{u_3}x_{4}^{u_4}$. Then there exist vectors ${\bf v}_{i}$, $1 \leq i \leq 22$, in $\mathbb{N}^{4}$ such that for all $w>0$, the toric ideal $I({\bf n}+w{\bf v}_{i})$ is a complete intersection whenever the entries of ${\bf n}+w{\bf v}_{i}$ are relatively prime.

\end{theorem}

\noindent \textbf{Proof.} By \cite[Theorem 6]{Kraft} $n_{1}=a_{2}(a_{3}u_{4}+u_{3}a_{4})$, $n_{2}=a_{1}(a_{3}u_{4}+u_{3}a_{4})$, $n_{3}=a_{4}(a_{1}u_{2}+u_{1}a_{2})$, $n_{4}=a_{3}(a_{1}u_{2}+u_{1}a_{2})$. Let ${\bf v}_{1}=(a_{2}a_{3},a_{1}a_{3},a_{2}a_{4},a_{2}a_{3})$ and $B=\{f_{1},f_{2},f_{4}=x_{1}^{u_{1}+w}x_{2}^{u_{2}}-x_{3}^{u_{3}}x_{4}^{u_{4}+w}\}$. Then $n_{1}+wa_{2}a_{3}=a_{2}(a_{3}(u_{4}+w)+u_{3}a_{4})$, $n_{2}+wa_{1}a_{3}=a_{1}(a_{3}(u_{4}+w)+u_{3}a_{4})$, $n_{3}+wa_{2}a_{4}=a_{4}(a_{1}u_{2}+(u_{1}+w)a_{2})$ and $n_{4}+wa_{2}a_{3}=a_{3}(a_{1}u_{2}+(u_{1}+w)a_{2})$. By Proposition \ref{Complete1} for every $w>0$, the ideal $I({\bf n}+w{\bf v}_{1})$ is a complete intersection on $f_{1}$, $f_{2}$ and $f_{4}$ whenever ${\rm gcd}(n_{1}+wa_{2}a_{3},n_{2}+wa_{1}a_{3},n_{3}+wa_{2}a_{4},n_{4}+wa_{2}a_{3})=1$. Consider the vectors ${\bf v}_{2}=(a_{2}a_{3},a_{1}a_{3},a_{1}a_{4},a_{1}a_{3})$, ${\bf v}_{3}=(a_{2}a_{4},a_{1}a_{4},a_{2}a_{4},a_{2}a_{3})$, ${\bf v}_{4}=(a_{2}a_{4},a_{1}a_{4},a_{1}a_{4},a_{1}a_{3})$, ${\bf v}_{5}=(a_{2}(a_{3}+a_{4}),a_{1}(a_{3}+a_{4}),0,0)$ and ${\bf v}_{6}=(0,0,a_{4}(a_{1}+a_{2}),a_{3}(a_{1}+a_{2}))$. By Proposition \ref{Complete1} for every $w>0$, $I({\bf n}+w{\bf v}_{2})$ is a complete intersection on $f_{1}$, $f_{2}$ and $x_{1}^{u_{1}}x_{2}^{u_{2}+w}-x_{3}^{u_{3}}x_{4}^{u_{4}+w}$ whenever the entries of ${\bf n}+w{\bf v}_{2}$ are relatively prime, $I({\bf n}+w{\bf v}_{3})$ is a complete intersection on $f_{1}$, $f_{2}$ and $x_{1}^{u_{1}+w}x_{2}^{u_{2}}-x_{3}^{u_{3}+w}x_{4}^{u_{4}}$ whenever the entries of ${\bf n}+w{\bf v}_{3}$ are relatively prime, and $I({\bf n}+w{\bf v}_{4})$ is a complete intersection on $f_{1}$, $f_{2}$ and $x_{1}^{u_{1}}x_{2}^{u_{2}+w}-x_{3}^{u_{3}+w}x_{4}^{u_{4}}$ whenever the entries of ${\bf n}+w{\bf v}_{4}$ are relatively prime. Furthermore for all $w>0$, $I({\bf n}+w{\bf v}_{5})$ is a complete intersection on $f_{1}$, $f_{2}$ and $x_{1}^{u_{1}}x_{2}^{u_{2}}-x_{3}^{u_{3}+w}x_{4}^{u_{4}+w}$ whenever the entries of ${\bf n}+w{\bf v}_{5}$ are relatively prime, and $I({\bf n}+w{\bf v}_{6})$ is a complete intersection on $f_{1}$, $f_{2}$ and $x_{1}^{u_{1}+w}x_{2}^{u_{2}+w}-x_{3}^{u_{3}}x_{4}^{u_{4}}$ whenever the entries of ${\bf n}+w{\bf v}_{6}$ are relatively prime. Consider the vectors ${\bf v}_{7}=(a_{2}(a_{3}+a_{4}),a_{1}(a_{3}+a_{4}),a_{2}a_{4},a_{2}a_{3})$, ${\bf v}_{8}=(a_{2}(a_{3}+a_{4}),a_{1}(a_{3}+a_{4}),a_{4}(a_{1}+a_{2}),a_{3}(a_{1}+a_{2}))$, ${\bf v}_{9}=(0,0,a_{2}a_{4},a_{2}a_{3})$, ${\bf v}_{10}=(a_{2}a_{4},a_{1}a_{4},a_{4}(a_{1}+a_{2}),a_{3}(a_{1}+a_{2}))$, ${\bf v}_{11}=(a_{2}a_{3},a_{1}a_{3},a_{4}(a_{1}+a_{2}),a_{3}(a_{1}+a_{2}))$, ${\bf v}_{12}=(a_{2}(a_{3}+a_{4}),a_{1}(a_{3}+a_{4}),a_{1}a_{4},a_{1}a_{3})$, ${\bf v}_{13}=(0,0,a_{1}a_{4},a_{1}a_{3})$, ${\bf v}_{14}=(a_{2}a_{4},a_{1}a_{4},0,0)$ and ${\bf v}_{15}=(a_{2}a_{3},a_{1}a_{3},0,0)$. Using Proposition \ref{Complete1} we have that for all $w>0$, $I({\bf n}+w{\bf v}_{i})$, $7 \leq i \leq 15$, is a complete intersection whenever the entries of ${\bf n}+w{\bf v}_{i}$ are relatively prime. For instance $I({\bf n}+w{\bf v}_{9})$ is a complete intersection on the binomials $f_{1}$, $f_{2}$ and $x_{1}^{u_{1}+w}x_{2}^{u_2}-x_{3}^{u_{3}}x_{4}^{u_{4}}$. Consider the vectors ${\bf v}_{16}=(a_{3}u_{4}+u_{3}a_{4},a_{3}u_{4}+u_{3}a_{4},a_{4}(u_{1}+u_{2}),a_{3}(u_{1}+u_{2}))$, ${\bf v}_{17}=(0,a_{3}u_{4}+u_{3}a_{4},u_{2}a_{4},u_{2}a_{3})$, ${\bf v}_{18}=(a_{3}u_{4}+u_{3}a_{4},0,u_{1}a_{4},u_{1}a_{3})$, ${\bf v}_{19}=(a_{2}u_{4},a_{1}u_{4},0,a_{1}u_{2}+u_{1}a_{2})$, ${\bf v}_{20}=(a_{2}u_{3},a_{1}u_{3},a_{1}u_{2}+u_{1}a_{2},0)$, ${\bf v}_{21}=(a_{2}(a_{4}+u_{4}),a_{1}(a_{4}+u_{4}),0,a_{1}u_{2}+u_{1}a_{2})$ and ${\bf v}_{22}=(a_{2}(u_{3}+u_{4}),a_{1}(u_{3}+u_{4}),a_{1}u_{2}+u_{1}a_{2},a_{1}u_{2}+u_{1}a_{2})$. It is easy to see that for all $w>0$, the ideal $I({\bf n}+w{\bf v}_{i})$, $16 \leq i \leq 22$, is a complete intersection whenever the entries of ${\bf n}+w{\bf v}_{i}$ are relatively prime. For instance $I({\bf n}+w{\bf v}_{16})$ is a complete intersection on the binomials $f_{2}$, $f_{3}$ and $x_{1}^{a_{1}+w}-x_{2}^{a_{2}+w}$. \hfill $\square$

\begin{example} {\rm Let ${\bf n}=(93,124,195,117)$, then $I({\bf n})$ is a complete intersection on the binomials $x_{1}^{4}-x_{2}^{3}$, $x_{3}^{3}-x_{4}^{5}$ and $x_{1}^{9}x_{2}^{3}-x_{3}^{2}x_{4}^{7}$. Here $a_{1}=4$, $a_{2}=3$, $a_{3}=3$, $a_{4}=5$, $u_{1}=9$, $u_{2}=3$, $u_{3}=2$ and $u_{4}=7$. Consider the vector ${\bf v}_{1}=(9,12,15,9)$. For all $w \geq 0$ the ideal $I({\bf n}+w{\bf v}_{1})$ is a complete intersection on $x_{1}^{4}-x_{2}^{3}$, $x_{3}^{3}-x_{4}^{5}$ and $x_{1}^{9+w}x_{2}^{3}-x_{3}^{2}x_{4}^{w+7}$ whenever ${\rm gcd}(93+9w,124+12w,195+15w,117+9w)=1$. By Corollary \ref{Arithmetical} the monomial curve $C({\bf n}+w{\bf v}_{1})$ has Cohen-Macaulay tangent cone at the origin. Consider the vector ${\bf v}_{4}=(15,20,20,12)$ and the sequence ${\bf n}+9{\bf v}_{4}=(228,304,375,225)$. The toric ideal $I({\bf n}+9{\bf v}_{4})$ is a complete intersection on the binomials $x_{1}^{4}-x_{2}^{3}$, $x_{3}^{3}-x_{4}^{5}$ and $x_{1}^{21}x_{2}^{3}-x_{3}^{2}x_{4}^{22}$. Note that $x_{1}^{25}-x_{3}^{2}x_{4}^{22} \in I({\bf n}+9{\bf v}_{4})$, so $x_{3}^{2}x_{4}^{22} \in I({\bf n}+9{\bf v}_{4})_{*}$ and also $x_{3}^{2} \in I({\bf n}+9{\bf v}_{4})_{*}: \langle x_{4}\rangle$. If $C({\bf n}+9{\bf v}_{4})$ has Cohen-Macaulay tangent cone at the origin, then $x_{3}^{2} \in I({\bf n}+9{\bf v}_{4})_{*}$ a contradiction. Thus $C({\bf n}+9{\bf v}_{4})$ does not have a Cohen-Macaulay tangent cone at the origin.}
\end{example}

\begin{theorem} \label{Basic1} Let $I({\bf n})$ be a complete intersection ideal generated by the binomials $f_{1}=x_{1}^{a_1}-x_{2}^{a_2}$, $f_{2}=x_{3}^{a_3}-x_{4}^{a_4}$ and $f_{3}=x_{1}^{u_1}x_{2}^{u_2}-x_{3}^{u_3}x_{4}^{u_4}$. Consider the vector ${\bf d}=(a_{2}a_{3},a_{1}a_{3},a_{2}a_{4},a_{2}a_{3})$. If $C({\bf n})$ has Cohen-Macaulay tangent cone at the origin, then for every $w>0$ the monomial curve $C({\bf n}+w{\bf d})$ has Cohen-Macaulay tangent cone at the origin whenever the entries of ${\bf n}+w{\bf d}$ are relatively prime.

\end{theorem}

\noindent \textbf{Proof.}  Let $n_{1}={\rm min}\{n_{1},\ldots,n_{4}\}$ and also $a_{3}<a_{4}$. Without loss of generality we can assume that $u_{2} \leq a_{2}$ and $u_{3} \leq a_{3}$. By Theorem \ref{Basic-complete} for every $w>0$, the ideal $I({\bf n}+w{\bf d})$ is a complete intersection on $f_{1}$, $f_{2}$ and $f_{4}=x_{1}^{u_{1}+w}x_{2}^{u_2}-x_{3}^{u_3}x_{4}^{u_{4}+w}$ whenever the entries of ${\bf n}+w{\bf d}$ are relatively prime. Note that $n_{1}+wa_{2}a_{3}={\rm min}\{n_{1}+wa_{2}a_{3},n_{2}+wa_{1}a_{3},n_{3}+wa_{2}a_{4},n_{4}+wa_{2}a_{3}\}$. Suppose that $u_{3}>0$ and $u_{4}>0$. Assume that $u_{2}<a_{2}$. By Corollary \ref{Arithmetical} it holds that $a_{4}+u_{4} \leq u_{1}+u_{2}+a_{3}-u_{3}$ and therefore $$a_{4}+(u_{4}+w) \leq (u_{1}+w)+u_{2}+a_{3}-u_{3}.$$ Thus, from Corollary \ref{Arithmetical} again $C({\bf n}+w{\bf d})$ has Cohen-Macaulay tangent cone at the origin. Assume that $u_{2}=a_{2}$. Then, from Corollary \ref{Arithmetical}, we have that $a_{4}+u_{4} \leq u_{1}+a_{1}+a_{3}-u_{3}$ and therefore $a_{4}+(u_{4}+w) \leq (u_{1}+w)+a_{1}+a_{3}-u_{3}$. By Corollary \ref{Arithmetical} $C({\bf n}+w{\bf d})$ has Cohen-Macaulay tangent cone at the origin.

Suppose now that $u_{3}=0$. Then $\{f_{1}, f_{2}, f_{5}=x_{4}^{u_{4}+w}-x_{1}^{u_{1}+w}x_{2}^{u_{2}}\}$ is a minimal generating set of $I({\bf n}+w{\bf d})$. If $u_{2}=a_{2}$, then $\{f_{1}, f_{2}, x_{4}^{u_{4}+w}-x_{1}^{u_{1}+a_{1}+w}\}$ is a standard basis for $I({\bf n}+w{\bf d})$ with respect to the negative degree reverse lexicographical order with $x_{3}>x_{4}>x_{2}>x_{1}$. Thus $I({\bf n}+w{\bf d})_{*}$ is a complete intersection and therefore $C({\bf n}+w{\bf d})$ has Cohen-Macaulay tangent cone at the origin. Assume that $u_{2}<a_{2}$, then $u_{4} \leq u_{1}+u_{2}$ and therefore $u_{4}+w \leq (u_{1}+w)+u_{2}$. The set $\{f_{1}, f_{2}, f_{5}\}$ is a standard basis for $I({\bf n}+w{\bf d})$ with respect to the negative degree reverse lexicographical order with $x_{3}>x_{4}>x_{2}>x_{1}$. Thus $I({\bf n}+w{\bf d})_{*}$ is a complete intersection and therefore $C({\bf n}+w{\bf d})$ has Cohen-Macaulay tangent cone at the origin.\\ Suppose that $u_{4}=0$, so necessarily $u_{3}=a_{3}$. Then $\{f_{1}, f_{2}, x_{4}^{a_{4}+w}-x_{1}^{u_{1}+w}x_{2}^{u_{2}}\}$ is a minimal generating set of $I({\bf n}+w{\bf d})$. If $u_{2}=a_{2}$, then $\{f_{1}, f_{2}, x_{4}^{a_{4}+w}-x_{1}^{u_{1}+a_{1}+w}\}$ is a standard basis for $I({\bf n}+w{\bf d})$ with respect to the negative degree reverse lexicographical order with $x_{3}>x_{4}>x_{2}>x_{1}$. Thus $I({\bf n}+w{\bf d})_{*}$ is a complete intersection and therefore $C({\bf n}+w{\bf d})$ has Cohen-Macaulay tangent cone at the origin. Assume that $u_{2}<a_{2}$, then $a_{4} \leq u_{1}+u_{2}$ and therefore $a_{4}+w \leq (u_{1}+w)+u_{2}$. The set $\{f_{1}, f_{2}, x_{4}^{a_{4}+w}-x_{1}^{u_{1}+w}x_{2}^{u_{2}}\}$ is a standard basis for $I({\bf n}+w{\bf d})$ with respect to the negative degree reverse lexicographical order with $x_{3}>x_{4}>x_{2}>x_{1}$. Thus $I({\bf n}+w{\bf d})_{*}$ is a complete intersection and therefore $C({\bf n}+w{\bf d})$ has Cohen-Macaulay tangent cone at the origin. \hfill $\square$

\begin{theorem} \label{TangentAll} Let $I({\bf n})$ be a complete intersection ideal generated by the binomials $f_{1}=x_{1}^{a_1}-x_{2}^{a_2}$, $f_{2}=x_{3}^{a_3}-x_{4}^{a_4}$ and $f_{3}=x_{1}^{u_1}x_{2}^{u_2}-x_{3}^{u_3}x_{4}^{u_4}$. Consider the vectors ${\bf d}_{1}=(0,0,a_{2}a_{4},a_{2}a_{3})$ and ${\bf d}_{2}=(0,0,a_{1}a_{4},a_{1}a_{3})$. Then there exists a non-negative integer $w_{0}$ such that for all $w \geq w_{0}$, the monomial curves $C({\bf n}+w{\bf d}_{1})$ and $C({\bf n}+w{\bf d}_{2})$ have Cohen-Macaulay tangent cones at the origin whenever the entries of the corresponding sequence (${\bf n}+w{\bf d}_{1}$ for the first family and ${\bf n}+w{\bf d}_{2}$ for the second) are relatively prime.

\end{theorem}

\noindent \textbf{Proof.} Let $n_{1}={\rm min}\{n_{1},\ldots,n_{4}\}$ and $a_{3}<a_{4}$. Suppose that $u_{2} \leq a_{2}$ and $u_{3} \leq a_{3}$. By Theorem \ref{Basic-complete} for all $w \geq 0$, $I({\bf n}+w{\bf d}_{1})$ is a complete intersection on $f_{1}$, $f_{2}$ and $f_{4}=x_{1}^{u_{1}+w}x_{2}^{u_2}-x_{3}^{u_3}x_{4}^{u_4}$ whenever the entries of ${\bf n}+w{\bf d_{1}}$ are relatively prime. Remark that $n_{1}={\rm min}\{n_{1},n_{2},n_{3}+wa_{2}a_{4},n_{4}+wa_{2}a_{3}\}$. Let $w_{0}$ be the smallest non-negative integer greater than or equal to $u_{3}+u_{4}-u_{1}-u_{2}+a_{4}-a_{3}$. Then for every $w \geq w_{0}$ we have that $a_{4}+u_{4} \leq u_{1}+w+u_{2}+a_{3}-u_{3}$, so $u_{3}+u_{4}<u_{1}+w+u_{2}$. Let $G=\{f_{1}, f_{2}, f_{4}, f_{5}=x_{4}^{a_{4}+u_{4}}-x_{1}^{u_{1}+w}x_{2}^{u_{2}}x_{3}^{a_{3}-u_{3}}\}$. We will prove that for every $w \geq w_{0}$, $G$ is a standard basis for $I({\bf n}+w{\bf d}_{1})$ with respect to the negative degree reverse lexicographical order with $x_{3}>x_{4}>x_{2}>x_{1}$. Note that ${\rm LM}(f_{1})=x_{2}^{a_2}$, ${\rm LM}(f_{2})=x_{3}^{a_{3}}$, ${\rm LM}(f_{4})=x_{3}^{u_{3}}x_{4}^{u_{4}}$ and ${\rm LM}(f_{5})=x_{4}^{a_{4}+u_{4}}$. Therefore ${\rm NF}({\rm spoly}(f_{i},f_{j})|G) = 0$ as ${\rm LM}(f_{i})$ and ${\rm LM}(f_j)$ are relatively prime, for $(i,j) \in \{(1,2),(1,4),(1,5),(2,5)\}$. We compute ${\rm spoly}(f_{2},f_{4})=-f_{5}$, so ${\rm NF}({\rm spoly}(f_{2},f_{4})|G)=0$. Next we compute ${\rm spoly}(f_{4},f_{5})=x_{1}^{u_{1}+w}x_{2}^{u_{2}}x_{3}^{a_{3}}-x_{1}^{u_{1}+w}x_{2}^{u_{2}}x_{4}^{a_{4}}$. Then ${\rm LM}({\rm spoly}(f_{4},f_{5}))=x_{1}^{u_{1}+w}x_{2}^{u_{2}}x_{3}^{a_{3}}$ and only ${\rm LM}(f_{2})$ divides ${\rm LM}({\rm spoly}(f_{4},f_{5}))$. Also ${\rm ecart}({\rm spoly}(f_{4},f_{5}))=a_{4}-a_{3}={\rm ecart}(f_{2})$. Then ${\rm spoly}(f_{2},{\rm spoly}(f_{4},f_{5}))=0$ and ${\rm NF}({\rm spoly}(f_{4},f_{5})|G)=0$. Thus the minimal number of generators for $I({\bf n}+w{\bf d}_{1})_{*}$ is either three or four, so from \cite[Theorem 3.1]{Shi} for every $w \geq w_{0}$, $C({\bf n}+w{\bf d}_{1})$ has Cohen-Macaulay tangent cone at the origin whenever the entries of ${\bf n}+w{\bf d}_{1}$ are relatively prime. 

By Theorem \ref{Basic-complete} for all $w \geq 0$, $I({\bf n}+w{\bf d}_{2})$ is a complete intersection on $f_{1}$, $f_{2}$ and $f_{6}=x_{1}^{u_{1}}x_{2}^{u_{2}+w}-x_{3}^{u_3}x_{4}^{u_4}$ whenever the entries of ${\bf n}+w{\bf d_{2}}$ are relatively prime. Remark that $n_{1}={\rm min}\{n_{1},n_{2},n_{3}+wa_{1}a_{4},n_{4}+wa_{1}a_{3}\}$. For every $w \geq w_{0}$ the set $H=\{f_{1}, f_{2}, f_{6}, x_{4}^{a_{4}+u_{4}}-x_{1}^{u_{1}}x_{2}^{u_{2}+w}x_{3}^{a_{3}-u_{3}}\}$ is a standard basis for $I({\bf n}+w{\bf d}_{2})$ with respect to the negative degree reverse lexicographical order with $x_{3}>x_{4}>x_{2}>x_{1}$. Thus the minimal number of generators for $I({\bf n}+w{\bf d}_{2})_{*}$ is either three or four, so from \cite[Theorem 3.1]{Shi} for every $w \geq w_{0}$, $C({\bf n}+w{\bf d}_{2})$ has Cohen-Macaulay tangent cone at the origin whenever the entries of ${\bf n}+w{\bf d}_{2}$ are relatively prime. \hfill $\square$

\begin{example} {\rm Let ${\bf n}=(15,25,24,16)$, then $I({\bf n})$ is a complete intersection on the binomials $x_{1}^{5}-x_{2}^{3}$, $x_{3}^{2}-x_{4}^{3}$ and $x_{1}x_{2}-x_{3}x_{4}$. Here $a_{1}=5$, $a_{2}=3$, $a_{3}=2$, $a_{4}=3$, $u_{i}=1$, $1 \leq i \leq 4$. Note that $x_{4}^{4}-x_{1}x_{2}x_{3} \in I({\bf n})$, so, from Corollary \ref{Arithmetical}, $C({\bf n})$ does not have a Cohen-Macaulay tangent cone at the origin. Consider the vector ${\bf d}_{1}=(0,0,9,6)$. For every $w>0$ the ideal $I({\bf n}+w{\bf d}_{1})$ is a complete intersection on the binomials $x_{1}^{5}-x_{2}^{3}$, $x_{3}^{2}-x_{4}^{3}$ and $x_{1}^{w+1}x_{2}-x_{3}x_{4}$ whenever ${\rm gcd}(15,25,24+9w,16+6w)=1$. By Theorem \ref{TangentAll} for every $w \geq 1$, the monomial curve $C({\bf n}+w{\bf d}_{1})$ has Cohen-Macaulay tangent cone at the origin whenever ${\rm gcd}(15,25,24+9w,16+6w)=1$.}
\end{example}

The next example gives a family of complete intersection monomial curves supporting M. Rossi's problem, although their tangent cones are not Cohen-Macaulay. To prove it we will use the following proposition. 
\begin{proposition} \cite[Proposition 2.2]{bayer} \label{hilbert} Let $I \subset K[x_{1}, x_{2}, \ldots, x_{d}]$ be a monomial ideal and $I=\langle J, {\bf x}^{\bf u} \rangle$ for a monomial ideal $J$ and a monomial ${\bf x}^{\bf u}$. Let $p(I)$ denote the numerator $g(t)$ of the Hilbert Series for $K[x_{1}, x_{2}, \ldots, x_{d}]/I$. Then $p(I)=p(J)-t^{{\rm deg}({\bf x}^{\bf u})}p(J: \langle{\bf x}^{\bf u} \rangle)$.
\end{proposition}

\begin{example} {\rm Consider the family
$n_{1}=8m^2+6$, $n_{2}=20m^2+15$, $n_{3}=12m^2+15$ and $n_{4}=8m^2+10$, where $m \geq 1$ is an integer. The toric ideal $I({\bf n})$ is minimally generated by the binomials
$$x_{1}^{5}-x_{2}^{2}, x_{3}^{2}-x_{4}^{3}, x_{1}^{2m^{2}}x_{2}-x_{3}x_{4}^{2m^{2}}.$$ Consider the vector ${\bf v}_{1}=(4,10,6,4)$ and the family $n_{1}'=n_{1}+4w$, $n_{2}'=n_{2}+10w$, $n_{3}'=n_{3}+6w$, $n_{4}'=n_{4}+4w$ where $w \geq 0$ is an integer. Let ${\bf n}'=(n_{1}',n_{2}',n_{3}',n_{4}')$, then for all $w \geq 0$ the toric ideal $I({\bf n}')$ is minimally generated by the binomials
$$x_{1}^{5}-x_{2}^{2}, x_{3}^{2}-x_{4}^{3}, x_{1}^{2m^{2}+w}x_{2}-x_{3}x_{4}^{2m^{2}+w}$$ whenever ${\rm gcd}(n_{1}',n_{2}',n_{3}',n_{4}')=1$. Let $C_{m}({\bf n}')$ be the corresponding monomial curve. By Corollary \ref{Arithmetical} for all $w \geq 0$, the monomial curve $C_{m}({\bf n}')$ does not have Cohen-Macaulay tangent cone at the origin whenever ${\rm gcd}(n_{1}',n_{2}',n_{3}',n_{4}')=1$. We will show that for every $w \geq 0$, the Hilbert function of the ring $K[[t^{n_{1}'},\ldots,t^{n_{4}'}]]$ is non-decreasing whenever ${\rm gcd}(n_{1}',n_{2}',n_{3}',n_{4}')=1$. It suffices to prove that for every $w \geq 0$, the Hilbert function of $K[x_{1}, x_{2}, x_3, x_{4}]/I({\bf n}')_*$ is non-decreasing whenever ${\rm gcd}(n_{1}',n_{2}',n_{3}',n_{4}')=1$. The set $$G=\{x_{1}^{5}-x_{2}^2, x_{3}^{2}-x_{4}^{3}, x_{1}^{2m^{2}+w}x_{2}-x_{3}x_{4}^{2m^{2}+w}, x_{4}^{2m^{2}+w+3}-x_{1}^{2m^{2}+w}x_{2}x_{3},$$ $$x_{1}^{2m^{2}+w+5}x_{3}-x_{2}x_{4}^{2m^{2}+w+3}, x_{1}^{4m^{2}+2w+5}-x_{4}^{4m^{2}+2w+3}\}$$ is a standard basis for $I({\bf n}')$ with respect to the negative degree reverse lexicographical order with $x_{4}>x_{3}>x_{2}>x_{1}$. Thus $I({\bf n}')_*$ is generated by the set
$$\{x_{2}^2, x_{3}^{2}, x_{4}^{4m^{2}+2w+3}, x_{1}^{2m^{2}+w}x_{2}x_{3}, x_{1}^{2m^{2}+w}x_{2}-x_{3}x_{4}^{2m^{2}+w}, x_{2}x_{4}^{2m^{2}+w+3}\}.$$ Also $\langle {\rm LT}(I({\bf n}')_*)\rangle$ with respect to the aforementioned order can be written as, $$\langle {\rm LT}(I({\bf n}')_*)\rangle=\langle x_{2}^{2}, x_{3}^{2}, x_{4}^{4m^{2}+2w+3}, x_{2}x_{4}^{2m^{2}+w+3}, x_{3}x_{4}^{2m^{2}+w}, x_{1}^{2m^{2}+w}x_{2}x_{3}\rangle.$$ Since the Hilbert function of $K[x_{1}, x_{2}, x_3, x_{4}]/I({\bf n}')_*$ is equal to the Hilbert function of $K[x_{1}, x_{2}, x_3, x_{4}]/\langle {\rm LT}(I({\bf n}')_*)\rangle$, it is sufficient to compute the Hilbert function of the latter.
Let $$J_{0}=\langle {\rm LT}(I({\bf n}')_*)\rangle, J_{1}=\langle x_{2}^{2}, x_{3}^{2}, x_{4}^{4m^{2}+2w+3}, x_{2}x_{4}^{2m^{2}+w+3}, x_{3}x_{4}^{2m^{2}+w} \rangle,$$ $$J_{2}=\langle x_{2}^{2}, x_{3}^{2}, x_{4}^{4m^{2}+2w+3}, x_{2}x_{4}^{2m^{2}+w+3} \rangle, J_{3}=\langle x_{2}^{2}, x_{3}^{2}, x_{4}^{4m^{2}+2w+3} \rangle.$$ Remark that $J_i= \langle J_{i+1}, q_i \rangle$, where $q_0=x_1^{2m^{2}+w}x_{2}x_{3}$, $q_1=x_{3}x_{4}^{2m^{2}+w}$ and $q_2=x_{2}x_{4}^{2m^{2}+w+3}$. We apply Proposition \ref{hilbert} to the ideal $J_i$ for $0 \leq i \leq 2$, so
\begin{equation} \label{recursion}
p(J_i)=p(J_{i+1})-t^{{\rm deg}(q_{i})}p(J_{i+1}:\langle q_i \rangle).
\end{equation}
Note that ${\rm deg}(q_{0})=2m^{2}+w+2$, ${\rm deg}(q_{1})=2m^{2}+w+1$ and ${\rm deg}(q_{2})=2m^{2}+w+4$.
In this case, it holds that
$J_{1}:\langle q_{0} \rangle=\langle x_{2}, x_{3},x_{4}^{2m^{2}+w} \rangle$,
$J_{2}: \langle q_{1} \rangle=\langle x_{2}^{2}, x_{3}, x_{4}^{2m^{2}+w+3}, x_{2}x_{4}^3 \rangle$ and $J_{3}: \langle q_{2} \rangle= \langle x_{2}, x_{3}^2, x_{4}^{2m^{2}+w} \rangle$. We have that
$$p(J_3)=(1-t)^{3}(1+3t+4t^{2}+\cdots+4t^{4m^{2}+2w+2}+3t^{4m^{2}+2w+3}+t^{4m^{2}+2w+4}).$$ Substituting all these recursively in Equation (\ref{recursion}),
we obtain that the Hilbert series of
$K[x_{1}, x_{2}, x_3, x_{4}]/J_0$ is
$$\frac{1+3t+4t^{2}+\cdots+4t^{2m^{2}+w}+3t^{2m^{2}+w+1}+t^{2m^{2}+w+2}+t^{2m^{2}+w+3}+t^{4m^{2}+2w+2}}{1-t}.$$
Since the numerator does not have any negative coefficients, the Hilbert function of $K[x_{1}, x_{2}, x_3, x_{4}]/J_{0}$ is non-decreasing whenever ${\rm gcd}(n_{1}',n_{2}',n_{3}',n_{4}')=1$.}
\end{example}

\section{The case (B)}

In this section we assume that after permuting variables, if necessary, $S=\{x_{1}^{a_1}-x_{2}^{a_2}, x_{3}^{a_{3}}-x_{1}^{u_1}x_{2}^{u_2}, x_{4}^{a_{4}}-x_{1}^{v_1}x_{2}^{v_2}x_{3}^{v_3}\}$ is a minimal generating set of $I({\bf n})$. Proposition \ref{Complete2} will be useful in the proof of Theorem \ref{Basic-complete3}.
\begin{proposition} \label{Complete2} Let $B=\{f_{1}=x_{1}^{b_1}-x_{2}^{b_2}, f_{2}=x_{3}^{b_{3}}-x_{1}^{c_1}x_{2}^{c_2}, f_{3}=x_{4}^{b_{4}}-x_{1}^{m_1}x_{2}^{m_2}x_{3}^{m_3}\}$ be a set of binomials in $K[x_{1},\ldots,x_{4}]$, where $b_{i} \geq 1$ for all $1 \leq i \leq 4$, at least one of $c_{1}$, $c_{2}$ is non-zero and at least one of $m_{1}$, $m_{2}$ and $m_{3}$ is non-zero.
Let $n_{1}=b_{2}b_{3}b_{4}$, $n_{2}=b_{1}b_{3}b_{4}$, $n_{3}=b_{4}(b_{1}c_{2}+c_{1}b_{2})$, $n_{4}=m_{3}(b_{1}c_{2}+b_{2}c_{1})+b_{3}(b_{1}m_{2}+m_{1}b_{2})$. If ${\rm gcd}(n_{1},\ldots,n_{4})=1$, then $I({\bf n})$ is a complete intersection ideal generated by the binomials $f_{1}$, $f_{2}$, $f_{3}$.

\end{proposition}

\noindent \textbf{Proof.} Consider the vectors ${\bf d}_{1}=(b_{1},-b_{2},0,0)$, ${\bf d}_{2}=(-c_{1},-c_{2},b_{3},0)$ and ${\bf d}_{3}=(-m_{1},-m_{2},-m_{3},b_{4})$. Clearly ${\bf d}_{i} \in {\rm ker}_{\mathbb{Z}}(n_{1},\ldots,n_{4})$ for $1 \leq i \leq 3$, so the lattice $L=\sum_{i=1}^{3} \mathbb{Z} {\bf d}_{i}$ is a subset of ${\rm ker}_{\mathbb{Z}}(n_{1},\ldots,n_{4})$. Let $$M=\begin{pmatrix}
 b_{1} & -c_{1} & -m_{1} \\
 -b_{2} & -c_{2} & -m_{2} \\
 0 & b_{3} & -m_{3} \\
 0 & 0 & b_{4}
\end{pmatrix},$$ then the rank of $M$ equals $3$. We will prove that the invariant factors $\delta_{1}$, $\delta_{2}$ and $\delta_{3}$ of $M$ are all equal to $1$. The greatest common divisor of all non-zero $3 \times 3$ minors of $M$ equals the greatest common divisor of the integers $n_{1}$, $n_{2}$, $n_{3}$ and $n_{4}$. But ${\rm gcd}(n_{1},\ldots,n_{4})=1$, so $\delta_{1} \delta_{2} \delta_{3}=1$ and therefore $\delta_{1}=\delta_{2}=\delta_{3}=1$. Note that the rank of the lattice ${\rm ker}_{\mathbb{Z}}(n_{1},\ldots,n_{4})$ is 3 and also equals the rank of $L$. By \cite[Lemma 8.2.5]{Vil} we have that $L={\rm ker}_{\mathbb{Z}}(n_{1},\ldots,n_{4})$. Now the transpose $M^{t}$ of $M$ is mixed dominating. By \cite[Theorem 2.9]{FS} the ideal $I({\bf n})$ is a complete intersection on $f_{1}$, $f_{2}$ and $f_{3}$.
 \hfill $\square$

\begin{theorem} \label{Basic-complete3} Let $I({\bf n})$ be a complete intersection ideal generated by the binomials $f_{1}=x_{1}^{a_1}-x_{2}^{a_2}$, $f_{2}=x_{3}^{a_{3}}-x_{1}^{u_1}x_{2}^{u_2}$ and $f_{3}=x_{4}^{a_{4}}-x_{1}^{v_1}x_{2}^{v_2}x_{3}^{v_3}$. Then there exist vectors ${\bf b}_{i}$, $1 \leq i \leq 22$, in $\mathbb{N}^{4}$ such that for all $w >0$, the toric ideal $I({\bf n}+w{\bf b}_{i})$ is a complete intersection whenever the entries of ${\bf n}+w{\bf b}_{i}$ are relatively prime.

\end{theorem}

\noindent \textbf{Proof.} By \cite[Theorem 6]{Kraft} $n_{1}=a_{2}a_{3}a_{4}$, $n_{2}=a_{1}a_{3}a_{4}$, $n_{3}=a_{4}(a_{1}u_{2}+u_{1}a_{2})$, $n_{4}=v_{3}(a_{1}u_{2}+a_{2}u_{1})+a_{3}(a_{1}v_{2}+v_{1}a_{2})$. Let ${\bf b}_{1}=(a_{2}a_{3},a_{1}a_{3},a_{1}u_{2}+u_{1}a_{2},a_{2}a_{3})$ and consider the set $B=\{f_{1},f_{2},f_{4}=x_{4}^{a_{4}+w}-x_{1}^{v_{1}+w}x_{2}^{v_{2}}x_{3}^{v_{3}}\}$. Then $n_{1}+wa_{2}a_{3}=a_{2}a_{3}(a_{4}+w)$, $n_{2}+wa_{1}a_{3}=a_{1}a_{3}(a_{4}+w)$, $n_{3}+w(a_{1}u_{2}+u_{1}a_{2})=(a_{4}+w)(a_{1}u_{2}+u_{1}a_{2})$ and $n_{4}+wa_{2}a_{3}=v_{3}(a_{1}u_{2}+a_{2}u_{1})+a_{3}(a_{1}v_{2}+(v_{1}+w)a_{2})$. By Proposition \ref{Complete2} for every $w>0$, the ideal $I({\bf n}+w{\bf b}_{1})$ is a complete intersection on $f_{1}$, $f_{2}$ and $f_{4}$ whenever the entries of ${\bf n}+w{\bf b}_{1}$ are relatively prime. Consider the vectors ${\bf b}_{2}=(a_{2}a_{3},a_{1}a_{3},a_{1}u_{2}+u_{1}a_{2},a_{1}a_{3})$, ${\bf b}_{3}=(a_{2}a_{3},a_{1}a_{3},a_{1}u_{2}+u_{1}a_{2},a_{1}u_{2}+u_{1}a_{2})$, ${\bf b}_{4}=(0,0,0,a_{3}(a_{1}+a_{2}))$, ${\bf b}_{5}=(0,0,0,a_{1}u_{2}+a_{2}u_{1}+a_{2}a_{3})$ and ${\bf b}_{6}=(0,0,0,a_{1}u_{2}+a_{2}u_{1}+a_{1}a_{3})$. By Proposition \ref{Complete2} for every $w>0$, $I({\bf n}+w{\bf b}_{2})$ is a complete intersection on $f_{1}$, $f_{2}$ and $x_{4}^{a_{4}+w}-x_{1}^{v_{1}}x_{2}^{v_{2}+w}x_{3}^{v_3}$ whenever the entries of ${\bf n}+w{\bf b}_{2}$ are relatively prime, $I({\bf n}+w{\bf b}_{3})$ is a complete intersection on $f_{1}$, $f_{2}$ and $x_{4}^{a_{4}+w}-x_{1}^{v_{1}}x_{2}^{v_{2}}x_{3}^{v_{3}+w}$ whenever the entries of ${\bf n}+w{\bf b}_{3}$ are relatively prime, and $I({\bf n}+w{\bf b}_{4})$ is a complete intersection on $f_{1}$, $f_{2}$ and $x_{4}^{a_{4}}-x_{1}^{v_{1}+w}x_{2}^{v_{2}+w}x_{3}^{v_{3}}$ whenever the entries of ${\bf n}+w{\bf b}_{4}$ are relatively prime. Furthermore for every $w>0$, $I({\bf n}+w{\bf b}_{5})$ is a complete intersection on $f_{1}$, $f_{2}$ and $x_{4}^{a_{4}}-x_{1}^{v_{1}+w}x_{2}^{v_{2}}x_{3}^{v_{3}+w}$ whenever the entries of ${\bf n}+w{\bf b}_{5}$ are relatively prime, and $I({\bf n}+w{\bf b}_{6})$ is a complete intersection on $f_{1}$, $f_{2}$ and $x_{4}^{a_{4}}-x_{1}^{v_{1}}x_{2}^{v_{2}+w}x_{3}^{v_{3}+w}$ whenever the entries of ${\bf n}+w{\bf b}_{6}$ are relatively prime. Consider the vectors ${\bf b}_{7}=(a_{2}a_{3},a_{1}a_{3},a_{1}u_{2}+u_{1}a_{2},a_{3}(a_{1}+a_{2}))$, ${\bf b}_{8}=(a_{2}a_{3},a_{1}a_{3},a_{1}u_{2}+u_{1}a_{2},a_{1}u_{2}+u_{1}a_{2}+a_{2}a_{3})$, ${\bf b}_{9}=(a_{2}a_{3},a_{1}a_{3},a_{1}u_{2}+u_{1}a_{2},a_{1}u_{2}+u_{1}a_{2}+a_{1}a_{3})$, ${\bf b}_{10}=(0,0,0,a_{1}u_{2}+a_{2}u_{1}+a_{3}(a_{1}+a_{2}))$, ${\bf b}_{11}=(a_{2}a_{3},a_{1}a_{3},a_{1}u_{2}+u_{1}a_{2},0)$, ${\bf b}_{12}=(0,0,0,a_{2}a_{3})$, ${\bf b}_{13}=(0,0,0,a_{1}a_{3})$, ${\bf b}_{14}=(0,0,0,a_{1}u_{2}+a_{2}u_{1})$ and ${\bf b}_{15}=(a_{2}a_{3},a_{1}a_{3},a_{1}u_{2}+u_{1}a_{2},a_{1}u_{2}+u_{1}a_{2}+a_{3}(a_{1}+a_{2}))$. Using Proposition \ref{Complete2} we have that for all $w>0$, the ideal $I({\bf n}+w{\bf b}_{i})$, $7 \leq i \leq 15$, is a complete intersection whenever the entries of ${\bf n}+w{\bf b}_{i}$ are relatively prime. Finally consider the vectors ${\bf b}_{16}=(a_{3}a_{4},a_{3}a_{4},a_{4}(u_{1}+u_{2}),v_{3}(u_{1}+u_{2})+a_{3}(v_{1}+v_{2}))$, ${\bf b}_{17}=(0,a_{3}a_{4},a_{4}u_{2},u_{2}v_{3}+a_{3}v_{2})$, ${\bf b}_{18}=(a_{3}a_{4},0,a_{4}u_{1},u_{1}v_{3}+v_{1}a_{3})$, ${\bf b}_{19}=(a_{2}a_{4},a_{1}a_{4},a_{2}a_{4},a_{2}v_{3}+a_{1}v_{2}+v_{1}a_{2})$, ${\bf b}_{20}=(a_{2}a_{4},a_{1}a_{4},a_{1}a_{4},a_{1}v_{3}+a_{1}v_{2}+v_{1}a_{2})$, ${\bf b}_{21}=(a_{2}a_{4},a_{1}a_{4},a_{4}(a_{1}+a_{2}),v_{3}(a_{1}+a_{2})+a_{1}v_{2}+v_{1}a_{2})$ and ${\bf b}_{22}=(0,0,a_{4}(a_{1}+a_{2}),v_{3}(a_{1}+a_{2})+a_{3}(a_{1}+a_{2}))$. It is easy to see that for all $w>0$, the ideal $I({\bf n}+w{\bf b}_{i})$, $16 \leq i \leq 22$, is a complete intersection whenever the entries of ${\bf n}+w{\bf b}_{i}$ are relatively prime. For instance $I({\bf n}+w{\bf b}_{22})$ is a complete intersection on the binomials $f_{1}$, $x_{3}^{a_{3}}-x_{1}^{u_{1}+w}x_{2}^{u_{2}+w}$ and $x_{4}^{a_{4}}-x_{1}^{v_{1}+w}x_{2}^{v_{2}+w}x_{3}^{v_{3}}$.
\hfill $\square$

\begin{example} {\rm Let ${\bf n}=(231,770,1023,674)$, then $I({\bf n})$ is a complete intersection on the binomials $x_{1}^{10}-x_{2}^{3}$, $x_{3}^{7}-x_{1}^{11}x_{2}^{6}$ and $x_{4}^{11}-x_{1}x_{2}^{8}x_{3}$. Here $a_{1}=10$, $a_{2}=3$, $a_{3}=7$, $a_{4}=11$, $u_{1}=11$, $u_{2}=6$, $v_{1}=1$, $v_{2}=8$ and $v_{3}=1$. Consider the vector ${\bf b}_{22}=(0,0,143,104)$, then for all $w \geq 0$ the ideal $I({\bf n}+w{\bf b}_{22})$ is a complete intersection on $x_{1}^{10}-x_{2}^{3}$, $x_{3}^{7}-x_{1}^{11+w}x_{2}^{6+w}$ and $x_{4}^{11}-x_{1}^{1+w}x_{2}^{8+w}x_{3}$ whenever ${\rm gcd}(231,770,1023+143w,674+104w)=1$. In fact, $I({\bf n}+w{\bf b}_{22})$ is minimally generated by $x_{1}^{10}-x_{2}^{3}$, $x_{3}^{7}-x_{1}^{11+w}x_{2}^{6+w}$ and $x_{4}^{11}-x_{1}^{11+w}x_{2}^{5+w}x_{3}$. Remark that $231={\rm min}\{231,770,1023+143w,674+104w\}$. The set $\{x_{1}^{10}-x_{2}^{3},x_{3}^{7}-x_{1}^{11+w}x_{2}^{6+w},x_{4}^{11}-x_{1}^{11+w}x_{2}^{5+w}x_{3}\}$ is a standard basis for $I({\bf n}+w{\bf b}_{22})$ with respect to the negative degree reverse lexicographical order with $x_{4}>x_{3}>x_{2}>x_{1}$. So $I({\bf n}+w{\bf b}_{22})_{*}$ is a complete intersection on $x_{2}^{3}$, $x_{3}^{7}$ and $x_{4}^{11}$, and therefore for every $w \geq 0$ the monomial curve $C({\bf n}+w{\bf b}_{22})$ has Cohen-Macaulay tangent cone at the origin whenever ${\rm gcd}(231,770,1023+143w,674+104w)=1$. Let ${\bf b}_{16}=(77,77,187,80)$. For every $w \geq 0$, $I({\bf n}+w{\bf b}_{16})$ is a complete intersection on $x_{1}^{10+w}-x_{2}^{3+w}$, $x_{3}^{7}-x_{1}^{11}x_{2}^{6}$ and $x_{4}^{11}-x_{1}x_{2}^{8}x_{3}$ whenever ${\rm gcd}(231+77w,770+77w,1023+187w,674+80w)=1$. Note that $231+77w={\rm min}\{231+77w,770+77w,1023+187w,674+80w\}$. For $0 \leq w \leq 5$ the set $\{x_{1}^{10+w}-x_{2}^{3+w},x_{3}^{7}-x_{1}^{11}x_{2}^{6},x_{4}^{11}-x_{1}^{11+w}x_{2}^{5-w}x_{3}\}$ is a standard basis for $I({\bf n}+w{\bf b}_{16})$ with respect to the negative degree reverse lexicographical order with $x_{4}>x_{3}>x_{2}>x_{1}$. Thus $I({\bf n}+w{\bf b}_{16})_{*}$ is minimally generated by $\{x_{2}^{3+w},x_{3}^{7},x_{4}^{11}\}$, so for $0 \leq w \leq 5$ the monomial curve $C({\bf n}+w{\bf b}_{16})$ has Cohen-Macaulay tangent cone at the origin whenever ${\rm gcd}(231+77w,770+77w,1023+187w,674+80w)=1$. Suppose that there is $w \geq 6$ such that $C({\bf n}+w{\bf b}_{16})$ has Cohen-Macaulay tangent cone at the origin. Then $x_{2}^{8}x_{3} \in I({\bf n}+w{\bf b}_{16})_{*}: \langle x_{1} \rangle$ and therefore $x_{2}^{8}x_{3} \in I({\bf n}+w{\bf b}_{16})_{*}$. Thus $x_{2}^{8}x_{3}$ is divided by $x_{2}^{3+w}$, a contradiction. Consequently for every $w \geq 6$ the monomial curve $C({\bf n}+w{\bf b}_{16})$ does not have Cohen-Macaulay tangent cone at the origin whenever ${\rm gcd}(231+77w,770+77w,1023+187w,674+80w)=1$.}
\end{example}

\begin{theorem} \label{theorcomp} Let $I({\bf n})$ be a complete intersection ideal generated by the binomials $f_{1}=x_{1}^{a_1}-x_{2}^{a_2}$, $f_{2}=x_{3}^{a_{3}}-x_{1}^{u_1}x_{2}^{u_2}$ and $f_{3}=x_{4}^{a_{4}}-x_{1}^{v_1}x_{2}^{v_2}x_{3}^{v_{3}}$. Consider the vector ${\bf d}=(0,0,a_{4}(a_{1}+a_{2}),v_{3}(a_{1}+a_{2})+a_{3}(a_{1}+a_{2}))$. Then there exists a non-negative integer $w_{1}$ such that for all $w \geq w_{1}$, the ideal $I({\bf n}+w{\bf d})_{*}$ is a complete intersection whenever the entries of ${\bf n}+w{\bf d}$ are relatively prime.

\end{theorem}

\noindent \textbf{Proof.} By Theorem \ref{Basic-complete3} for all $w \geq 0$, the ideal $I({\bf n}+w{\bf d})$ is minimally generated by $G=\{f_{1},f_{4}=x_{3}^{a_{3}}-x_{1}^{u_{1}+w}x_{2}^{u_{2}+w},f_{5}=x_{4}^{a_{4}}-x_{1}^{v_{1}+w}x_{2}^{v_{2}+w}x_{3}^{v_{3}}\}$ whenever the entries of ${\bf n}+w{\bf d}$ are relatively prime. Let $w_{1}$ be the smallest non-negative integer greater than or equal to ${\rm max}\{\frac{a_{3}-u_{1}-u_{2}}{2},\frac{a_{4}-v_{1}-v_{2}-v_{3}}{2}\}$. Then $a_{3} \leq u_{1}+u_{2}+2w_{1}$ and $a_{4} \leq v_{1}+v_{2}+v_{3}+2w_{1}$. It is easy to prove that for every $w \geq w_{1}$ the set $G$ is a standard basis for $I({\bf n}+w{\bf d})$ with respect to the negative degree reverse lexicographical order with $x_{4}>x_{3}>x_{2}>x_{1}$. Note that ${\rm LM}(f_{1})$ is either $x_{1}^{a_{1}}$ or $x_{2}^{a_{2}}$, ${\rm LM}(f_{4})=x_{3}^{a_{3}}$ and ${\rm LM}(f_{5})=x_{4}^{a_{4}}$. By \cite[Lemma 5.5.11]{GP} $I({\bf n}+w{\bf d})_{*}$ is generated by the least homogeneous summands of the elements in the standard basis $G$. Thus for all $w \geq w_{1}$, the ideal $I({\bf n}+w{\bf d})_{*}$ is a complete intersection whenever the entries of ${\bf n}+w{\bf d}$ are relatively prime. \hfill $\square$

\begin{proposition} \label{AlmostShibuta} Let $I({\bf n})$ be a complete intersection ideal generated by the binomials $f_{1}=x_{1}^{a_1}-x_{2}^{a_2}$, $f_{2}=x_{3}^{a_{3}}-x_{1}^{u_1}x_{2}^{u_2}$ and $f_{3}=x_{4}^{a_{4}}-x_{1}^{v_1}x_{2}^{v_2}$, where $v_{1}>0$ and $v_{2}>0$. Assume that $a_{2}<a_{1}$, $a_{3}<u_{1}+u_{2}$, $v_{2}<a_{2}$ and $a_{1}+v_{1} \leq a_{2}-v_{2}+a_{4}$. Then there exists a vector ${\bf b}$ in $\mathbb{N}^{4}$ such that for all $w \geq 0$, the ideal $I({\bf n}+w{\bf b})_{*}$ is almost complete intersection whenever the entries of ${\bf n}+w{\bf b}$ are relatively prime.

\end{proposition}

\noindent \textbf{Proof.} From the assumptions we deduce that $v_{1}+v_{2}<a_{4}$. Consider the vector ${\bf b}=(a_{2}a_{3},a_{1}a_{3},a_{1}u_{2}+u_{1}a_{2},a_{2}a_{3})$. For every $w \geq 0$ the ideal $I({\bf n}+w{\bf b})$ is a complete intersection on $f_{1}$, $f_{2}$ and $f_{4}=x_{4}^{a_{4}+w}-x_{1}^{v_{1}+w}x_{2}^{v_{2}}$ whenever the entries of ${\bf n}+w{\bf b}$ are relatively prime. We claim that the set $G=\{f_{1},f_{2},f_{4},f_{5}=x_{1}^{a_{1}+v_{1}+w}-x_{2}^{a_{2}-v_{2}}x_{4}^{a_{4}+w}\}$ is a standard basis for $I({\bf n}+w{\bf b})$ with respect to the negative degree reverse lexicographical order with $x_{3}>x_{2}>x_{1}>x_{4}$. Note that ${\rm LM}(f_{1})=x_{2}^{a_2}$, ${\rm LM}(f_{2})=x_{3}^{a_{3}}$, ${\rm LM}(f_{4})=x_{1}^{v_{1}+w}x_{2}^{v_{2}}$ and ${\rm LM}(f_{5})=x_{1}^{a_{1}+v_{1}+w}$. Also ${\rm spoly}(f_{1},f_{4})=-f_{5}$. It suffices to show that ${\rm NF}({\rm spoly}(f_{4},f_{5})|G)=0$. We compute ${\rm spoly}(f_{4},f_{5})=x_{2}^{a_2}x_{4}^{a_{4}+w}-x_{1}^{a_1}x_{4}^{a_{4}+w}$. Then ${\rm LM}({\rm spoly}(f_{4},f_{5}))=x_{2}^{a_2}x_{4}^{a_{4}+w}$ and only ${\rm LM}(f_1)$ divides ${\rm LM}({\rm spoly}(f_{4},f_{5}))$. Moreover ${\rm ecart}({\rm spoly}(f_{4},f_{5}))=a_{1}-a_{2}={\rm ecart}(f_{1})$. So ${\rm spoly}(f_{1},{\rm spoly}(f_{4},f_{5}))=0$ and also ${\rm NF}({\rm spoly}(f_{4},f_{5})|G)=0$. Thus \begin{enumerate}
\item If $a_{1}+v_{1}<a_{2}-v_{2}+a_{4}$, then $I({\bf n}+w{\bf b})_{*}$ is minimally generated by $\{x_{2}^{a_{2}},x_{3}^{a_{3}},x_{1}^{v_{1}+w}x_{2}^{v_2},x_{1}^{a_{1}+v_{1}+w}\}$.
\item If $a_{1}+v_{1}=a_{2}-v_{2}+a_{4}$, then $I({\bf n}+w{\bf b})_{*}$ is minimally generated by $\{x_{2}^{a_{2}},x_{3}^{a_{3}},x_{1}^{v_{1}+w}x_{2}^{v_2},f_{5}\}$. \hfill $\square$

\end{enumerate}


\end{document}